\documentclass[12pt]{article}
\usepackage{epsfig}
\newcommand{\bfx}{\mathbf{x}}

\begin{document}

\title{\bf Scaling universalities of $k$th-nearest neighbor
           distances on closed manifolds}
\author{
Allon G.~Percus\thanks{E-mail address: percus@lanl.gov} \\
{\normalsize CIC--3 and Center for Nonlinear Studies, MS-B258}\\
{\normalsize Los Alamos National Laboratory}\\
{\normalsize Los Alamos, NM 87545, USA}\\
\\
and \\
\\
Olivier C.~Martin\thanks{E-mail address: martino@ipno.in2p3.fr} \\
{\normalsize Division de Physique Th\'eorique\thanks{Unit\'e de Recherche
des Universit\'es Paris XI et Paris VI associ\'ee au C.N.R.S.}}\\
{\normalsize Institut de Physique Nucl\'eaire, Universit\'e Paris-Sud}\\
{\normalsize F--91406 Orsay Cedex, France}\\
}

\date{}
\maketitle
\begin{abstract}
\medskip
Take $N$ sites distributed randomly and uniformly on a smooth closed surface.
We express
the expected distance $\langle D_k(N)\rangle$ from an arbitrary point on the
surface to its $k$th-nearest neighboring site, in terms of the function
$A(l)$ giving the
area of a disc of radius $l$ about that point.  We then find two
universalities.  First, for a flat surface, where $A(l)=\pi l^2$, the
$k$-dependence and the $N$-dependence separate in $\langle D_k(N)\rangle$.
All $k$th-nearest neighbor distances thus have the same scaling law in
$N$.
Second, for a curved surface, the average $\int\langle D_k(N)\rangle\,d\mu$
over the surface is a topological invariant at leading and subleading order
in a large $N$ expansion.  The $1/N$ scaling series then depends, up
through $O(1/N)$, only on the surface's
topology and not on its precise shape.  We discuss the case
of higher dimensions ($d>2$), and also interpret our results using Regge
calculus.

\end{abstract}

\bigskip
\par
\centerline{Key words: stochastic geometry, nearest neighbors}
\vspace{1cm}
\centerline{Submitted to {\it Advances in Applied Mathematics\/}, February 1998}

\newpage
\baselineskip=20pt

\section{Introduction}
\label{sect_introduction}
Many problems arising in applied mathematics involve the distance between
neighboring sites in a space.  One frequently wishes to calculate
distances to a nearest neighbor, second-nearest neighbor, and
more generally, $k$th-nearest neighbor.  Examples in computational
geometry and optimization abound, ranging from random packing of
spheres to minimum spanning trees.  Such problems also occur naturally in 
physics, and have been considered in applications ranging from stellar
dynamics \cite{Chandrasekhar}, to interactions in liquid systems,
to cellular objects such as
foams and random lattices \cite{ItzyksonDrouffe}.

Here we consider the case of $N$ sites placed randomly, with a
uniform distribution, on a 2-D surface of fixed area.  Let the random
variable $D_k(N)$ represent the distance between a given point $\bfx$ and
its $k$th-nearest site.  The expectation value $\langle
D_k(N)\rangle$ taken over the ensemble of randomly placed
sites --- and in fact all moments $\langle D^\alpha_k(N)\rangle$ ---
then exhibit some surprising properties.  
First of all,
when the surface is flat $\langle D_k(N)\rangle$ may be written, up
to corrections exponentially small in $N$, as
\begin{displaymath}
\langle D_k(N)\rangle \approx \frac{1}{\sqrt{\pi}}\,
\frac{(k-1/2)!}{(k-1)!}\,\frac{N!}{(N+1/2)!}\mbox{.}
\end{displaymath}
This means that the $k$-dependence and the $N$-dependence separate; the
large $N$ scaling law for $\langle D_k(N)\rangle$ is independent of $k$.
Geometrically, the meaning of this universality is far from obvious.
Furthermore, the property is not restricted to two dimensions, and turns
out to be equally valid for flat spaces of any dimension.
Second of all, when the surface is curved, while this $k$-independence no
longer holds, one finds on the other hand another universality: if $\langle
D_k(N)\rangle$ is written in terms of a $1/N$ series, and averaged
over the entire surface, the leading and subleading coefficients of the
expansion give topological invariants.  Thus to $O(1/N)$, the large
$N$ scaling law depends not on the detailed
shape of the surface but only on the surface's genus.

In this paper we explore these universalities.
We start by expressing $\langle D_k(N)\rangle$
in terms of the area, $A(l)$, of a disc of radius $l$ on an arbitrary surface.
We observe that in the special case where $A(l)$ consists only of a
power of $l$,
the scaling law for $\langle D_k(N)\rangle$ exhibits the universality
in $k$, i.e., is independent of $k$.  Then, we give a relation between
$A(l)$ and the Gaussian curvature of the surface, and find the leading
correction terms in the $1/N$ power series for $\langle D_k(N)\rangle$
as functions of the curvature.  This results in the topological invariance
at $O(1/N)$.  We discuss higher order terms as well, and the case of higher
dimensions.  Finally, we show that a Regge calculus
approach provides a simple means of obtaining this topological invariance
for the case of polyhedral (non-smooth) surfaces.

\section{Preliminaries}
\label{sect_prelim}

\begin{figure}[!b]
\begin{center}
\epsfig{file=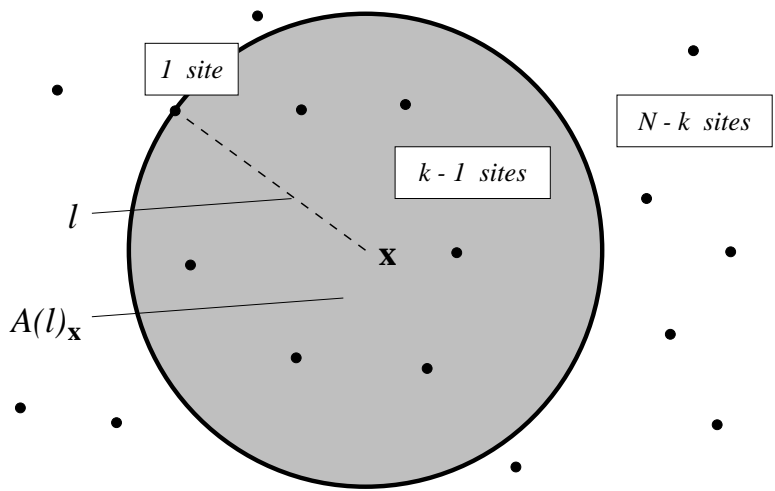,width=4in}
\end{center}
\caption{A point $\bfx$ and its neighborhood within a distance $l$.
Area $A(l)_\bfx$ is the probability of finding a site within distance $l$
of $\bfx$.
Here, $k-1$ sites lie within this region, the $k$th-nearest site lies
at distance $l$ exactly, and the remaining $N-k$ sites lie outside the
region.}
\label{fig_neighbors}
\end{figure}

Take any point $\bfx$, and consider $P[D_k(N)=l]_\bfx$, the probability
density that the point's $k$th-nearest neighboring site lies at a distance
$l$ from it.  This is equal to the probability density of having $k-1$
(out of $N$) sites {\it within\/} distance $l$, one site (out of $N-k+1$)
{\it at\/} distance $l$, and the remaining $N-k$ sites
{\it beyond\/} distance $l$.  Let us choose units so that our surface has
total area 1.  Since sites are distributed uniformly over the surface,
the probability of a site lying within distance $l$ is then simply
the area $A(l)_\bfx$ of a disc of radius $l$ about point $\bfx$ on the
surface.  This is shown in Figure \ref{fig_neighbors}.
Dropping the argument $\bfx$ (in order to simplify the notation),
we may then write
\begin{displaymath}
P[D_k(N)=l] = {N \choose k-1}\,\left[ A(l)\right]^{k-1}\times\,
{N-k+1 \choose 1}\,\frac{dA(l)}{dl}\,\times\,\left[ 1-A(l)
\right]^{N-k},
\end{displaymath}
giving the expectation value (first moment)
\begin{eqnarray*}
\langle D_k(N)\rangle &=& \int_0^\infty P[D_k(N)=l]\,l\,dl\\
&=& \frac{N!}{(N-k)!\,(k-1)!}\,\int_0^\infty l\,\left[ A(l)\right]^{k-1}
\left[ 1-A(l)\right]^{N-k}\,\frac{dA(l)}{dl}\,dl.
\end{eqnarray*}
Using the variable transformation $w=A(l)$, this may be written in
terms of the inverse function $A^{-1}(w)$ as
\begin{displaymath}
\langle D_k(N)\rangle = \frac{N!}{(N-k)!\,(k-1)!}\,\int_0^1
A^{-1}(w)\,w^{k-1} (1-w)^{N-k}\, dw.
\end{displaymath}

If $A^{-1}(w)$ admits the power series expansion in $w$:
\begin{equation}
\label{eq_series1}
A^{-1}(w) = w^\gamma\,\sum_{j=0}^\infty c_j w^j,\qquad
\mbox{for some }\gamma\in[0,1),
\end{equation}
then
\begin{equation}
\label{eq_integral}
\langle D_k(N)\rangle = \frac{N!}{(N-k)!\,(k-1)!}\,\sum_{j=0}^\infty c_j
\int_0^1 w^{k+j+\gamma-1} (1-w)^{N-k}\, dw.
\end{equation}
Recognizing the integral as the Beta function ${\rm B}(k+j+\gamma,N-k+1)=
(k+j+\gamma-1)!\,(N-k)!/$ $(N+j+\gamma)!$,
\begin{equation}
\label{eq_series2}
\langle D_k(N)\rangle = \sum_{j=0}^\infty c_j\,
\frac{(k+j+\gamma-1)!}{(k-1)!}\,\frac{N!}{(N+j+\gamma)!}.
\end{equation}

Several comments are in order concerning $\langle D_k(N)\rangle$.  First
of all, although we restrict ourselves to discussing the first moment
of $D_k(N)$, we could in fact consider any moment $\langle
D^\alpha_k(N)\rangle$ by taking $\left[ A^{-1}(w)\right]^\alpha$ instead
of $A^{-1}(w)$ in (\ref{eq_series1}).  Doing so would alter $\gamma$ and
the $c_j$'s, but would not change our results qualitatively.  Second of all,
there is no loss of generality in taking our total surface area
to be unity; scaling this area by a constant (or even, as might be more
intuitive in statistical physics, by $N$) would provide
only a trivial scaling factor in our results.  Third of all, we could
imagine that the point $\bfx$ we consider
is itself an $(N+1)$th site.  This is simply a question of nomenclature:
the problem of finding the expected distance {\it from an arbitrary point
to its kth-nearest site\/}, for a system of $N$ sites, is equivalent
to the problem of finding the expected distance {\it between kth-nearest
neighboring sites\/}, for a system of $N+1$ sites.

We now turn to the properties of $A^{-1}(w)$, and their consequences on
$\langle D_k(N)\rangle$.

\section{Flat Surfaces}
\label{sect_flat}

On a flat surface, if we could neglect edge effects, the area included
within distance $l$ would simply be $A(l)=\pi l^2$.  In that case,
we would have $A^{-1}(w)=\sqrt{w/\pi}$, and so from (\ref{eq_series1}) and
(\ref{eq_series2}),
\begin{equation}
\label{eq_noexp}
\langle D_k(N)\rangle = \frac{1}{\sqrt{\pi}}\,\frac{(k-1/2)!}{(k-1)!}\,
\frac{N!}{(N+1/2)!}.
\end{equation}
There would thus be a complete separation of the $k$-dependence and the
$N$-dependence.

As we are working with a surface of fixed (unit) area, however, we
cannot avoid considering edge or finite size effects.  Let us restrict ourselves
to the case where the surface is everywhere locally Euclidean within
some minimum neighborhood of radius $l_0>0$.  (The simplest example of this is
a unit square with periodic boundary conditions, for which $l_0=1/2$.
Clearly, many other constructions are possible.)
Any required modification to the $A^{-1}(w)$ expression in
(\ref{eq_series1}) then concerns only
$w$ greater than $w_0\equiv A(l_0)$.  Correspondingly, (\ref{eq_integral})
remains valid, up to remainder terms from the region of integration
$w_0\le w\le 1$.  Since the
$(1-w)^{N-k-1}$ term in the integral is bounded above by $(1-w_0)^{N-k-1}$
within this region, these remainder terms are exponentially small in $N$.
Equation (\ref{eq_noexp}) is thus still correct {\it to all
orders\/} in a $1/N$ series expansion, and may be written as
\begin{displaymath}
\langle D_k(N)\rangle = \frac{1}{\sqrt{\pi}}\,\frac{(k-1/2)!}{(k-1)!}\,
\frac{1}{\sqrt{N}}\left[ 1 - \frac{3}{8N} +
O\left(\frac{1}{N^2}\right)\right],
\end{displaymath}
where all orders in the series are independent of $k$.  We therefore
see that the large $N$ scaling law for
$k$th-nearest neighbor distances on a 2-D flat surface without a boundary
exhibits the universality in $k$ to all orders in $1/N$.

The same holds true for flat manifolds in any dimension $d$.  We assume
there is some $l_0$ such that the
volume included within distance $l<l_0$ is simply
the volume of a $d$-dimensional ball:
\begin{eqnarray*}
A(l)&=&\frac{\pi^{d/2}\,l^d}{(d/2)!}, \qquad\mbox{or}\\
A^{-1}(w)&=&\frac{1}{\sqrt{\pi}}\left[ w\left( \frac{d}{2}\right) !
\right]^{1/d}.
\end{eqnarray*}
As before, this condition allow us to write
$\langle D_k(N)\rangle$ up to remainder terms that are exponentially
small in $N$, so from (\ref{eq_series2}),
\begin{eqnarray*}
\langle D_k(N)\rangle &\approx&
\frac{\left[ (d/2)!\right]^{1/d}}{\sqrt{\pi}}\,
\frac{(k-1+1/d)!}{(k-1)!}\,\frac{N!}{(N+1/d)!}\\
&=&
\frac{\left[ (d/2)!\right]^{1/d}}{\sqrt{\pi}}\,
\frac{(k-1+1/d)!}{(k-1)!}\, N^{-1/d}\left[
1- \frac{1/d + 1/d^2}{2N} + O\left( \frac{1}{N^2}\right)\right].
\end{eqnarray*}
Thus for flat spaces without a boundary, of any dimension $d$,
the universality in $k$ holds to all orders in $1/N$: the $k$-dependence
and the $N$-dependence separate.

\begin{figure}[t]
\begin{center}
\begin{tabular}{c@{\hspace{0.2in}}c}
\epsfig{file=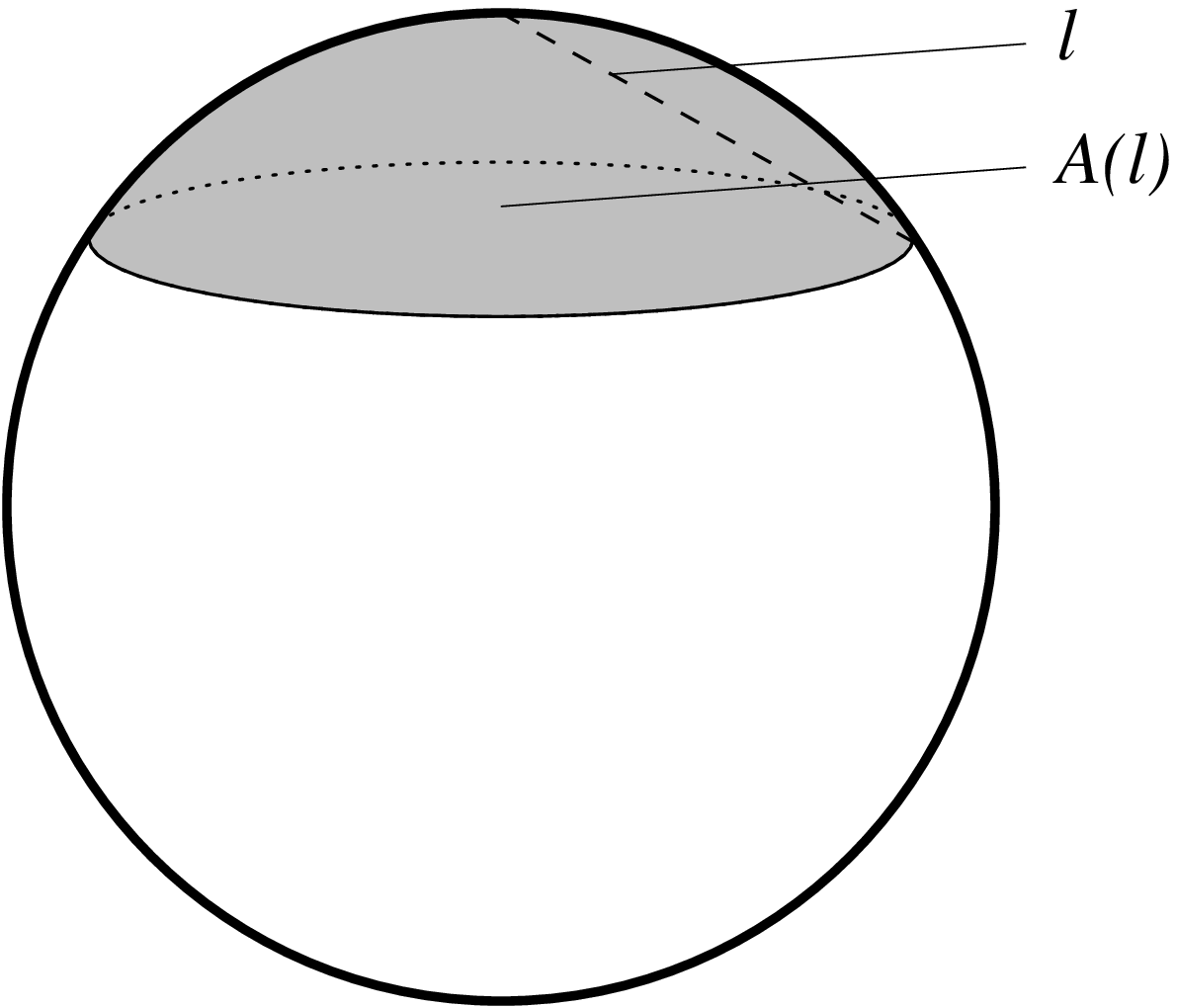,width=2.5in} &
\epsfig{file=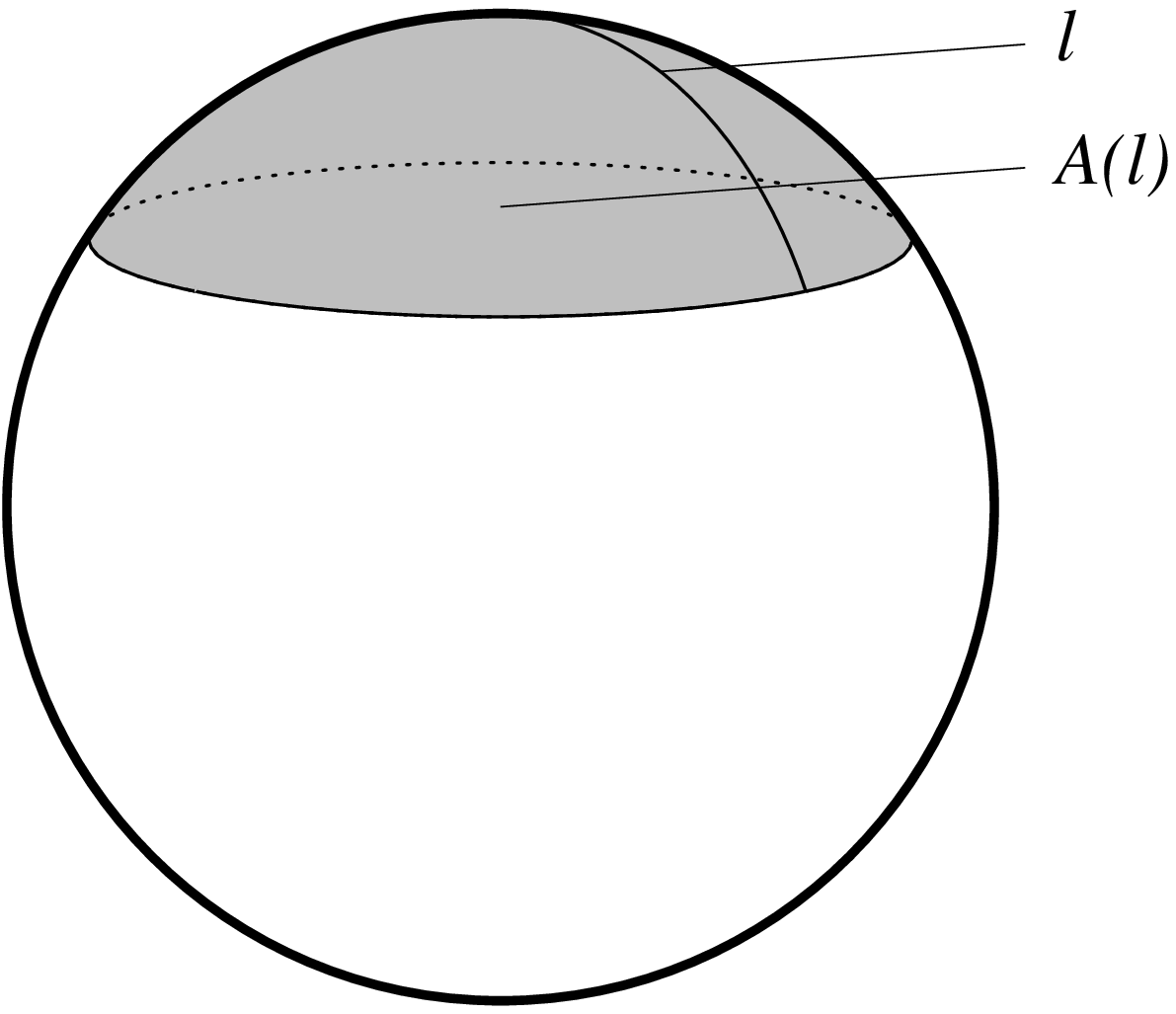,width=2.5in} \\
(a) $l$ represents chord length distance\hspace{0.35in} &
(b) $l$ represents arc length distance\hspace{0.35in} \\
\end{tabular}
\end{center}
\caption{A spherical surface, where $A(l)$ denotes the area of a
spherical cap spanned by $l$.  In (a), $l$ represents the chord length
distance; in (b), $l$ represents the arc length distance.}
\label{fig_sphere}
\end{figure}

It may be interesting to consider a slight variation on the problem,
giving this universality {\it exactly\/} and not only to all
orders.  Take the case of a spherical surface embedded in
3-D Euclidean space, with the usual measure of area
over the sphere, but with a peculiar
sort of ``distance'': rather than the conventional choice of the
arc length (geodesic) metric, use the chord length.  (See Figure
\ref{fig_sphere}a.)  For
a chord of length $l$ originating at a pole of the sphere, the area
of the spherical cap spanned by it is simply $A(l)=\pi l^2$.
The $k$th-nearest neighbor distance properties using chord length
``distance'' on this curved surface then appear analogous to those
on a flat surface.  
There is, however, one important distinction.  The relevant
threshold $w_0$ for edge effects is in this case $w_0=\pi (2R)^2$,
where $R$ is the radius of the sphere.  Since $\pi (2R)^2$ is exactly
equal to the total surface area of the sphere, it is set to 1.  Equation
(\ref{eq_integral}) thus requires no corrections at all, and so the
universality in (\ref{eq_noexp}) is exact.

\section{Curved Surfaces}
\label{sect_curved}

Now consider the case of a surface with intrinsic curvature, with the
distance defined in terms of a metric, i.e., along
geodesics of the surface.  Let us begin with a spherical surface,
this time letting $l$ represent arc length; the area of the spherical cap
spanned by an arc originating at a pole of the sphere (see Figure
\ref{fig_sphere}b) is given by
\begin{eqnarray*}
A(l)_{\rm sphere} &=& 2\pi R^2 \left[ 1-\cos\frac{l}{R}\right] \\
&=& 4\pi R^2 \sin^2\frac{l}{2R}.
\end{eqnarray*}
If the total surface area, $4\pi R^2$, is normalized to 1,
\begin{eqnarray}
A(l)_{\rm sphere} &=& \sin^2\sqrt{\pi}l,\qquad\mbox{so}\nonumber\\
A^{-1}(w)_{\rm sphere} &=& \frac{\sin^{-1}\sqrt{w}}{\sqrt{\pi}}\nonumber\\
&\simeq&\sqrt{\frac{w}{\pi}}\,\sum_{j=0}^\infty\,\frac{w^j}{2j+1}\,
\frac{(2j)!}{2^{2j}(j!)^2}.
\label{eq_area_sphere}
\end{eqnarray}
As in the case of the chord length ``distance'', this $A^{-1}(w)$
expression
is exact everywhere for $0\le w\le 1$.  Equations (\ref{eq_series1}) and
(\ref{eq_series2}) then require no corrections, and we find
\begin{eqnarray}
\langle D_k(N)\rangle_{\rm sphere} &=&
\frac{1}{\sqrt{\pi}}\,
\sum_{j=0}^\infty\,\frac{1}{2j+1}\,\frac{(2j)!}{2^{2j}(j!)^2}\,
\frac{(k+j-1/2)!}{(k-1)!}\,\frac{\sqrt{N}\,N!}{(N+j+1/2)!}\nonumber\\
&=&
\frac{1}{\sqrt{\pi}}\,\frac{(k-1/2)!}{(k-1)!}\,\frac{1}{\sqrt{N}}\left[
1+ \frac{4k-7}{24N} + O\left( \frac{1}{N^2}\right)\right].
\label{eq_nn_sphere}
\end{eqnarray}
Clearly, the $k$ universality does not apply here: the $O(1/N)$ coefficient
explicitly contains $k$.

Another sort of universality, however, is found when we turn to the more
general case of an arbitrary closed surface, i.e., an abstract 2-D
manifold with no boundary.  Given a
smooth surface,
we may introduce a system of curvilinear coordinates $u$ and $v$ and
write (at least piecewise on the surface) the differential length element
$ds$ in the conformal, orthogonal form \cite{Kreyszig}:
\begin{equation}
\label{eq_ds}
ds^2 = f(u,v)\,[du^2 + dv^2].
\end{equation}
The Gaussian curvature $K(u,v)$ of the surface is then expressed
in terms of the function $f(u,v)$ by
\begin{equation}
\label{eq_gausscurv}
K = \frac{1}{2f^3}\left[\left(\frac{\partial f}{\partial u}\right)^2
+ \left(\frac{\partial f}{\partial v}\right)^2
- f\frac{\partial^2 f}{\partial u^2}
- f\frac{\partial^2 f}{\partial v^2}\right].
\end{equation}
What is $A(l)$ on this surface? To find out, we first determine 
the manifold's geodesic lines.  For
$ds$ given by (\ref{eq_ds}), we may use the geodesic equation \cite{Kreyszig}:
\begin{equation}
\label{eq_geodesic}
\frac{d^2 u}{ds^2}+\frac{1}{2f}\,\frac{\partial f}{\partial u}
\left[\left(\frac{du}{ds}\right)^2 - \left(\frac{dv}{ds}\right)^2\right]
+ \frac{1}{f}\,\frac{\partial f}{\partial v}\,\frac{du}{ds}\,\frac{dv}{ds}
=0.
\end{equation}
Let us expand $u$ and $v$ as functions of distance $s$ from an initial
point, along a fixed geodesic:
\begin{eqnarray}
u(s) &=& u_0 + s u_0' + \frac{s^2}{2} u_0'' + \cdots\quad\mbox{and}\nonumber\\
v(s) &=& v_0 + s v_0' + \frac{s^2}{2} v_0'' + \cdots ,
\label{eq_expansion}
\end{eqnarray}
where $u_0\equiv u(0)$, $u_0'\equiv u'(0)$, etc., and likewise for $v$.
Then, expanding $f(u,v)$ in terms of $u$ and $v$ and substituting
(\ref{eq_expansion}),
\begin{eqnarray*}
f(u,v)
&=& f(u_0,v_0) + s \left[ u_0' f_u(u_0,v_0) + v_0' f_u(u_0,v_0)\right]
+ s^2 \left[\frac{u_0''}{2}f_u(u_0,v_0) + \frac{v_0''}{2}f_v(u_0,v_0)\right.\\
&& \qquad\left. + \frac{(u_0')^2}{2}f_{uu}(u_0,v_0) +
u_0' v_0' f_{uv}(u_0,v_0) +
\frac{(v_0')^2}{2}f_{vv}(u_0,v_0)\right] + O(s^3),
\end{eqnarray*}
where subscripts on $f$ denote partial derivatives.

Using (\ref{eq_ds}) and (\ref{eq_geodesic}), we can solve for all
but three of the coefficients in (\ref{eq_expansion}).  Let us choose
$u_0$, $v_0$ and $u_0'$ to be these three.  Now, consider the area $A(l)$
about the point $(u_0,v_0)$.  For $ds$ given in (\ref{eq_ds}), the
differential surface element will be $d\mu=f\,du\,dv$, so:
\begin{eqnarray}
A(l) &=& \int f\,du\,dv\nonumber\\
&=& \int f\,J\left(\frac{u,v}{s,u_0'}\right)ds\,du_0'\nonumber\\
&=& \int f\, \left|\frac{\partial u}{\partial s}
\frac{\partial v}{\partial u_0'} - \frac{\partial u}{\partial u_0'}
\frac{\partial v}{\partial s}\right|\,ds\,du_0'.
\label{eq_area_int}
\end{eqnarray}
The limits of integration over $s$ are $0$ and $l$;
the limits of integration over $u_0'$, which may be found from
(\ref{eq_ds}), are $-\sqrt{1/f(u_0,v_0)}$ and $\sqrt{1/f(u_0,v_0)}$.
In evaluating the Jacobian some care must be taken, as a sign ambiguity
allows two solutions for the coefficients in (\ref{eq_expansion}).
$A(l)$ will be the sum of (\ref{eq_area_int}) evaluated at each of the
two solutions, ultimately causing all odd powers of $l$ to vanish.

The result, after lengthy algebraic manipulations, may be written as
\begin{equation}
\label{eq_area_series1}
A(l)=\pi l^2\left[ 1-\frac{l^2}{24}\,\frac{1}{f^3}\,(f_u^2 + f_v^2
-f f_{uu} -f f_{vv}) + O(l^4)\right],
\end{equation}
where in order to avoid cluttering the notation, we have omitted the
$(u_0,v_0)$ arguments at which all functions are to be evaluated.
For the leading correction term in $A(l)$ given in
(\ref{eq_area_series1}), we recognize
the expression (\ref{eq_gausscurv}) for the Gaussian curvature $K$.
In retrospect, this is not surprising: by symmetry, the correction
series to $\pi l^2$ can contain only even powers of $l$,
and if we consider $A(l)$
as a geometric expansion about a flat space approximation, $K$ will be
the only scalar curvature quantity with dimensions of $l^{-2}$
\cite{Weinberg}.

With some perseverance,
one may carry the expansion in (\ref{eq_area_series1}) to higher orders,
obtaining
\begin{eqnarray}
A(l) &=& \pi l^2\left[ 1-\frac{l^2}{12}\,K + \frac{l^4}{720}\bigl(
2K^2-3\nabla^2K\bigr)\right.\nonumber\\
&& \qquad\left. - \frac{l^6}{161280}\bigl( 8K^3 -
3[10(\nabla K)^2 + 14K\nabla^2K - 5\nabla^4K]\bigr)
+ O(l^8)\right],
\label{eq_area_series2}
\end{eqnarray}
where $\nabla$ is the gradient operator.
We thus obtain a series expansion giving the area of a disc on 
a smooth 2-D
surface, in a form that depends only on intrinsic quantities, i.e., not
on the choice of coordinate system.

We may invert (\ref{eq_area_series2}) to obtain the power series
\begin{eqnarray}
\label{eq_area_inverse}
A^{-1}(w) &=& \sqrt{\frac{w}{\pi}}\,\Biggl[ 1+ \frac{K}{24\pi}\,w +
\frac{9K^2 + 4\nabla^2K}{1920\pi^2}\,w^2 +
\frac{15K^3 + 14K\nabla^2K - 2(\nabla K)^2 + \nabla^4K}
{21504\pi^3}\,w^3\Biggr.\nonumber\\
&& \qquad\Biggl. + O(w^4)\Biggr].
\end{eqnarray}
As an example, take the special case of a spherical surface, where the
Gaussian curvature is a constant $K=1/R^2$, or $K=4\pi$ for a unit
surface.  All derivatives of $K$ then vanish, leaving
\begin{displaymath}
A^{-1}(w)_{\rm sphere}=\sqrt{\frac{w}{\pi}}\left[ 1+ \frac{w}{6} +
\frac{3w^2}{40} + \frac{5w^3}{112} + O(w^4)\right],
\end{displaymath}
from which we recover the first few terms of our earlier result
(\ref{eq_area_sphere}).

Given an expression for $A^{-1}(w)$ on a general 2-D surface, we may now find
$\langle D_k(N)\rangle$ using (\ref{eq_series1}) and (\ref{eq_series2}).
It is helpful at this point to define the reduced variable
$\langle\tilde{D}_k(N)\rangle$ by dividing out the leading asymptotic
(large $N$) behavior from $\langle D_k(N)\rangle$.
Recalling the notation of (\ref{eq_series1}), $A^{-1}(w)=
w^\gamma\,\sum_j c_j w^j$ for some $\gamma\in[0,1)$, define:
\begin{eqnarray}
\langle\tilde{D}_k(N)\rangle &=& \langle D_k(N)\rangle\,\frac{1}{c_0}\,
\frac{(k-1)!}{(k+\gamma-1)!}\,N^\gamma\nonumber\\
&=& N^\gamma\,\sum_{j=0}^\infty\,\frac{c_j}{c_0}\,
\frac{(k+j+\gamma-1)!}{(k+\gamma-1)!}\,
\frac{N!}{(N+j+\gamma)!},
\label{eq_reduced}
\end{eqnarray}
so that $\lim_{N\to\infty}\langle\tilde{D}_k(N)\rangle = 1$ (this is
a consequence of Stirling's law).  $\langle\tilde{D}_k(N)\rangle$ then
corresponds to the $1/N$ power series that we have frequently seen appearing
in $\langle D_k(N)\rangle$, giving the corrections to its leading
(large $N$) behavior.  It is in this power series
that the interesting universality properties emerge.
Consider the average $\int\langle\tilde{D}_k(N)\rangle\,d\mu$ over
the entire surface, obtained --- via (\ref{eq_series1}) and
(\ref{eq_reduced}) --- from precisely the average of the
series coefficients in (\ref{eq_area_inverse}).  Examine, in particular,
the $O(w)$ term of (\ref{eq_area_inverse}).
By the Gauss-Bonnet theorem \cite{Kreyszig}, $\int K\,d\mu=2\pi\chi$ on any
closed surface, where $\chi$ is the Euler
characteristic of the surface, a topological invariant.
Up to leading corrections, then,
\begin{eqnarray}
\label{eq_chi}
\int A^{-1}(w)\,d\mu &=& \sqrt{\frac{w}{\pi}}\left[ 1+\frac{\chi}{12} w +
O(w^2)\right],\qquad\mbox{giving}\nonumber\\
\int\langle\tilde{D}_k(N)\rangle\,d\mu &=& \frac{\sqrt{N}\,N!}{(N+1/2)!}\,
\left[ 1+ \frac{\chi}{12}\,
\frac{k+1/2}{N+3/2} + O\left(\frac{1}{N^2}\right)\right]\nonumber\\
&=& 1+ \frac{\chi (2k+1)-9}{24N} + O\left( \frac{1}{N^2}\right).
\label{eq_nn_general}
\end{eqnarray}
We thus discover a different sort of universality from the one
we had in the case of flat space.  To $O(1/N)$, the scaling law
for $k$th-nearest neighbor distances depends only on the surface's topology,
and not on its detailed properties.

The Euler characteristic $\chi$ for a surface is related to its genus $g$
by $\chi=2(1-g)$.  Taking the torus as one example, $g=1$,
so $\chi=0$ and the
$k$-dependence in (\ref{eq_nn_general}) once again disappears, at least
at $O(1/N)$.  This is to be expected: a flat space with periodic boundary
conditions has, after all, the topology of a torus.  And conversely,
because of the topological invariance, all tori behave like flat space
to $O(1/N)$.  Taking the spherical surface as another example, $g=0$,
so $\chi=2$ and we
recover from (\ref{eq_nn_general}) the power series in (\ref{eq_nn_sphere}).

The properties of $\int\langle\tilde{D}_k(N)\rangle\,d\mu$ are far less clear
at higher orders in $1/N$.
Using the divergence theorem and integration by parts, we may obtain from
(\ref{eq_area_inverse})
\begin{equation}
\label{eq_area_inverse_av}
\int A^{-1}(w)\,d\mu=\int\sqrt{\frac{w}{\pi}}\left[ 1+ \frac{K}{24\pi}\,w +
\frac{3K^2}{640\pi^2}\,w^2 + \frac{15K^3 + 16K\nabla^2K}
{21504\pi^3}\,w^3 + O(w^4)\right]\,d\mu .
\end{equation}
If we looked only at terms up through $O(w^2)$, we might believe that this
series is simply, by analogy with (\ref{eq_area_sphere}), the expansion of
$\int (2/\sqrt{K}) \sin^{-1}\sqrt{Kw/4\pi}\,d\mu$.  Unfortunately,
starting at $O(w^3)$ we see this is not true, since the contributions
of curvature and its gradients
do not all vanish in the average over the surface!  Furthermore,
even for terms in (\ref{eq_area_inverse_av}) of the form $\int K^n\,d\mu$,
at $n>1$ there is no straightforward equivalent to the Gauss-Bonnet theorem;
the theorem is a direct consequence of the integrand's
linearity.  Thus for a general 2-D surface, a simplified form does not
appear to exist for the terms in $\int\langle\tilde{D}_k(N)\rangle\,d\mu$
beyond $O(1/N)$.  More particularly, the {\it only\/} case in which 
$\int\langle\tilde{D}_k(N)\rangle\,d\mu$ would be independent of $k$ 
beyond $O(1/N)$ is if the curvature
is identically equal to $0$, i.e., a flat surface.

Let us briefly consider the case of curved higher-dimensional manifolds.
The calculation is now far more complicated, as it is no longer possible
to write the metric tensor in a conformal form as we did in
(\ref{eq_ds}).  In addition, whereas in 2-D the only intrinsic
scalar quantity describing curvature is the Gaussian curvature $K$,
for $d>2$ there are $d(d-1)(d-2)(d+3)/12$ different such quantities
\cite{Weinberg}.  However, all of them except $K$ itself have dimensions
of higher order than $l^{-2}$.  It thus seems reasonable to conjecture that,
as we argued
in 2-D, the $O(l^2)$ correction term in $A(l)$ can only involve $K$.  
(Indeed, we have verified that this is true in 3-D.)  In that case, we
may rely on the example of the spherical surface --- easily generalized to
$d$ dimensions --- to provide us with the initial terms for a general
manifold:
\begin{eqnarray}
A(l) &=& \frac{\pi^{d/2}}{(d/2)!}\,l^d\,\left[ 1-\frac{d(d-1)}{d+2}\,
\frac{K}{6}\,l^2+O(l^4)\right],\qquad\mbox{or}\nonumber\\
A^{-1}(w) &=& \left[\frac{(d/2)!}{\pi^{d/2}}\right]^{1/d} w^{1/d}\,
\left[ 1+\frac{d-1}{d+2}\,\frac{K}{6}\,\left(\frac{(d/2)!}{\pi^{d/2}}
\right)^{2/d}\,w^{2/d}+O(w^{4/d})\right].
\label{eq_area_multid}
\end{eqnarray}
Note that $A^{-1}(w)$ now contains a series in $w^{2/d}$ rather than in
$w$.  Appropriately modifying (\ref{eq_series1}), it may then be shown
that $\langle\tilde{D}_k(N)\rangle$ is in general given by a series in
$1/N^{1/d}$ for odd $d$, and $1/N^{2/d}$ for even $d$.

Consider, finally, the average $\int A^{-1}(w)\,d\mu$ over the manifold.
The higher-dimensional generalization of the Gauss-Bonnet theorem
\cite{Santalo}
involves an integrand of $O(1/l^d)$, or $O(1/w)$.  The leading correction
term from (\ref{eq_area_multid}), $\int K\,d\mu$, therefore cannot be
simplified further for $d>2$; the only term that could possibly give
rise to a topological invariant is the coefficient at $O(l^d)$, or
$O(w)$.  If $d$ is odd, it is rather certain that no topological invariant
will exist in the series.  If $d$ is even, the $O(w)$ term will
first contribute to the $\int\langle\tilde{D}_k(N)\rangle\,d\mu$ series
at $O(1/N)$ --- as in 2-D, although at higher dimensions this will no
longer be the {\it leading\/} correction term.  While one cannot rule out the
possibility
of obtaining a topological invariant at $O(1/N)$, the $O(w)$
term in $A^{-1}(w)$ is in general a complicated one involving
many different curvature scalars, and so this is far from obvious.
We leave it as an open question.

\section{Regge Calculus}
\label{sect_regge}
We have remarked that from a physical point of view it is natural,
in the 2-D case, for the leading
corrections in $A(l)$ to contain only the
Gaussian curvature, as this represents the leading deviation from
planarity.  Consequently, only the
mean curvature --- or, using the Gauss-Bonnet theorem, the Euler
characteristic $\chi$ --- matters in the $O(1/N)$ term of 
$\int\langle\tilde{D}_k(N)\rangle\,d\mu$.  We have seen using
differential methods (geodesics) that this physical picture is
indeed correct.  These methods apply to a smooth closed surface.  For
polyhedral surfaces, which are not smooth, we may in fact obtain a
similar result
more easily, using the non-differential method of Regge
calculus.  Consider
a polyhedron with a number of vertices, edges, and faces.
Following the work of Regge \cite{Regge} and others since then
\cite{CheegerMullerSchrader}, we observe that the curvature is
concentrated at the vertices and is measured by a deficit angle: if
$\theta_i$ is the sum of the angles incident on vertex $i$, 
the deficit angle at that vertex is $\Delta_i = 2 \pi - \theta_i$.
It may then be shown that the Gauss-Bonnet theorem, on polyhedra,
reduces to Euler's relation
$2 \pi \chi = \sum_i \Delta_i$.

Let $P$ be a polyhedron with 
a fixed number of vertices, and consider
the problem of finding the large $N$ scaling series
$\int\langle\tilde{D}_k(N)\rangle\,d\mu$ on $P$.
As $N \to \infty$, corrections to the flat space value about a given
point $\bfx$ arise only when $\bfx$ is near one of the vertices,
because only in that case can curvature (i.e., the deficit angle) enter into
the local calculation of $A(l)$ about $\bfx$.
It is then sufficient to understand the corrections associated with 
one vertex at a time.
Consider the neighborhood of a vertex $i$.  $A(l)$ receives a
correction from the flat space value, and by a simple geometric construction,
one can see that this correction is exactly proportional to the
deficit angle $\Delta_i$.
Correspondingly, the leading correction
term both in $A^{-1}(w)$ and in the scaling series will
be proportional to $\Delta_i$ for small deficit angles.
Now  sum over all the vertices $i = 1,\dots,N$, 
assuming that all the deficit angles are small. We then
find that the $O(1/N)$
term in $\int\langle\tilde{D}_k(N)\rangle\,d\mu$
is proportional to $\chi$, and we recover the topological
invariant derived in the case of a smooth manifold.

A word of caution is necessary, however.  It is tempting at this point
to take the limit where $P$ becomes a smooth manifold, expecting to recover
(\ref{eq_chi}).  Unfortunately this will not work; a direct computation
shows that the limit does not commute with the limit $N\to\infty$ taken
above, and the coefficient thus obtained at $O(1/N)$ will
not be the correct one.

\section{Conclusions}
\label{sect_conclusions}

Given $N$ sites distributed randomly and uniformly on a surface with no
boundaries, we have considered the properties of mean distances
to neighboring sites.  When the surface is flat, we have seen that in
the expression for the mean $k$th-nearest site distance, the $k$-dependence
and $N$-dependence separate.  The scaling law in $N$ for mean $k$th-nearest
neighbor distances is thus independent of $k$.  This universality
applies equally well to higher moments of the distances,
and to Euclidean manifolds in dimensions greater than 2.
For surfaces with curvature, while this general property is no longer valid,
we have found that when the $k$th-nearest neighbor distance is written as
a large $N$ expansion, averaged over the surface, the leading correction
term in the series is a topological invariant.  The scaling series thus
depends, to $O(1/N)$, on the genus of the manifold but not on its other
properties.

Although we have considered these universalities only for the moments of
point-to-point distances, similar properties hold for higher order
simplices such as areas of triangles associated with nearby points.
The problem is thus a natural one to consider further in the context
of random triangulations, foams or other physical problems
\cite{ItzyksonDrouffe}
tightly connected to geometry.

\section*{Acknowledgments}
\label{sect_ack}
We are grateful to E.~Bogomolny, J.~Houdayer, and C.~Kenyon
for sharing with us their
valuable insights on this topic, and to O.~Bohigas for having introduced
us to the problem.  AGP wishes to acknowledge the hospitality of the
Division de Physique Th\'eorique, Institut de Physique Nucl\'eaire, Orsay,
where much of this work was carried out.  OCM acknowledges support from
the Institut Universitaire de France.

\end{document}